
\input amstex
\documentstyle{amsppt}
\magnification =\magstep1
\hsize 6.5truein
\vsize 9truein
%\addto\tenpoint{\normalbaselineskip=20pt\normalbaselines}
\expandafter\redefine\csname logo\string@\endcsname{}
\topmatter
\title\nofrills{$L$-series associated to symmetric functions mod $N$ with applications related to $\zeta(3), \zeta(5)$}
\endtitle
\leftheadtext{$L$-series for symmetric functions mod $N$}
\rightheadtext{David Spring}
\author David Spring
\endauthor
\address Department of Mathematics, Glendon College, York University, 2275 Bayview Avenue,
Toronto, Ontario, Canada, M4N 3M6.\endaddress

\email dspring\@glendon.yorku.ca\endemail

\abstract We develop a new theory of $L$-series based on replacing Dirichlet characters mod $N$
by symmetric functions mod $N$ in order to calculate explicitly the sums of infinite series more closely related
to $\zeta(2n+1)$, specifically $\zeta(3),\zeta(5)$. This generalizes the corresponding theory of sums of $L$-series associated to Dirichlet characters.\endabstract

\TagsOnRight
\endtopmatter

\document
\flushpar {\it Keywords}\,: $\zeta(2n+1)$, $L$-series, symmetric function, Dirichlet character

\flushpar Mathematics Subject Classification 2010: 11M36

\newpage
\centerline{\S 1. Introduction}\vskip.25cm
In this paper we develop a new theory of $L$-series obtained by replacing Dirichlet characters
by what we denote as symmetric functions mod $N$ in order to calculate explicitly the sums of infinite series more closely related
to $\zeta(2n+1)$, more specifically $\zeta(3),\zeta(5)$, thereby obtaining new results in  the extensive literature on zeta functions. Symmetric functions are more
flexible than Dirichlet characters in that the homomorphism property of a Dirichlet character is not required, thus increasing the scope of
associated $L$-series whose sums can still be calculated explicitly. Two examples taken from \S 3 Theorem 3.5, which are not defined by Dirichlet characters,
and which seem to be new to the literature, are as follows.

$$
\aligned
L(3,\chi_8)&=\left(1+\frac{1}{2^3}+\frac{1}{3^3}-\frac{1}{5^3}-\frac{1}{6^3}-\frac{1}{7^3}\right)\\
&+\left(\frac{1}{9^3}+\frac{1}{10^3}+\frac{1}{11^3}-
\frac{1}{13^3}-\frac{1}{14^3}-\frac{1}{15^3}\right)\\
&+ \dots =\frac{\pi^3}{256}(6\sqrt{2}+1).\endaligned
$$

$$
\aligned
L(3,\chi_{12})&=\left(1+\frac{1}{2^3}+\frac{1}{3^3}+\frac{1}{4^3} +\frac{1}{5^3}-\frac{1}{7^3}-\frac{1}{8^3}-\frac{1}{9^3}-\frac{1}{10^3}- \frac{1}{11^3}\right)\\
&+\left(\frac{1}{13^3}+\frac{1}{14^3}+\frac{1}{15^3}+\frac{1}{16^3}+\frac{1}{17^3}-\frac{1}{19^3}-\frac{1}{20^3}-\frac{1}{21^3}-\frac{1}{22^3} -\frac{1}{23^3}\right)\\
&+\dots =\frac{\pi^3}{7776}(20\sqrt{3}+261).
\endaligned
$$
The techniques employed in Theorem 3.5 enable one to calculate particular infinite series obtained by modifying $\zeta(2n+1)$ to have a sequence of
$k$ successive $+$ signs followed by a sequence of $k$ successive $-$ signs (separated by a term with coefficient 0), repeated periodically ad infinitum.
As $k\to \infty$ one recovers the series $\zeta(2n+1)$. However, a closed formula for $\zeta(2n+1)$ based on this approach remains unknown. Formally one proceeds as follows, beginning with symmetric functions.\vskip.25cm

{\bf Definition 1.1.} Let $r\ge 1$, $N\ge 2$. A {\it symmetric function mod $N$} is a function
$$ \chi (=\chi(N,r))\: \{1,2,\dots ,N-1\}\to \bold C$$
such that the following properties obtain:

\item{($S_1$):} $\chi(N-a)=(-1)^r\chi(a)$, $1\le a\le N-1$. In particular, $\chi(N/2)=0$ for $r$ odd and $N$ even.\vskip.25cm

 \item{($S_2$):} (periodicity) The function $\chi$ extends to a function (same notation) $\chi\: \bold Z\to \bold C$ such that $\chi(kN)=0$
 and $\chi(kN+a) =\chi(a)$ for all $k\in \bold Z$ and for all $a\in \{1,2, \dots ,N-1\}$. \vskip.25cm

\flushpar Associated to a symmetric function $\chi=\chi(N,r)$ is the $L$-series, $L(r,\chi)=\sum_{n=1}^\infty\frac{\chi(n)}{n^r}$.\vskip.25cm

\flushpar{\bf Remark 1.2.} Note that $\chi(-m)=(-1)^r\chi(m)$ for all $m\in \bold Z$. Indeed from $(S_2)$ if $m=kN+a$, $a\in \{1,\dots ,N-1\}$
then $\chi(m)=\chi(a)$; $-m=-kN-a=-(k+1)N+N-a$, hence from $(S_1)$, $(S_2)$, $\chi(-m)=\chi(N-a)=(-1)^r\chi(a)=(-1)^r\chi(m)$. Thus $\chi$ is an
even (odd) periodic function if $r$ is even (odd). In this paper, the case $r$ odd is the main interest.
\vskip.25cm

\flushpar{\bf Remark 1.3.} Recall [3, p.\,115], [2, p.\,82] that a Dirichlet character mod $N$, $N\ge 2$, is a homomorphism of groups
$f\: \bold Z_N^{\times}\to \bold C^{\times}$, ($R^{\times}$ is the multiplicative group of units of a commutative ring $R$ with identity).
Extending $f$ (same notation) to $\bold Z$ by property $(S_2)$ above, such that also $f(a)=0$ if $\gcd(a,N)\ge 2$, it follows that $f\: \bold Z\to \bold C$
is (totally) multiplicative: $f(ab)=f(a)f(b)$, for all $a,b\in\bold Z$. Note that $(f(-1))^2=f(1)=1$, hence $f(-1)=\pm 1$. Consequently a Dirichlet character also
satisfies property $(S_1)$ above, and hence is a special case of a symmetric function mod $N$: $f(N-a)=f(-a)=(-1)^rf(a)$ for all $a$, where $r$ is odd if $f(-1)=-1$,
respectively $r$ is even if $f(-1)=1$ (cf. also the Historical Note 2.4 in \S2).

\flushpar{\bf Remark 1.4.} In general a symmetric function $\chi$ mod $N$ is not a
Dirichlet character since its restriction to $\bold Z_N^{\times}$ is not required to be a homomorphism of groups. For example, from $(S_1)$ the function values $\chi(a)$
can be chosen arbitrarily,
$1\le a<N/2$, $N\ge 3$. More explicitly, the symmetric function $\chi_{2m}$ defined in (1.2) below is not a Dirichlet character for all $m\ge 3$. Employing
Remark 1.3, we see
that the symmetric functions mod $N$ generalize Dirichlet characters mod $N$ in the following sense: A Dirichlet character $f$ mod $N$ defines an infinite
sequence of symmetric functions $\chi_f(N,r)$, indexed by integers $r\ge 1$, such that $\chi_f(a)=f(a)$ for all $a\in\bold Z$, and (i) $r$ is odd if $f(-1)=-1$; (ii)
$r$ is even if $f(-1)=1$. \vskip.25cm

\flushpar The $L$-series associated to a Dirichlet character $f$ mod $N$ (viewed as a symmetric function $\chi_f(=\chi_f(N,r))$ mod $N$, as in Remark (1.4)
is exactly one of the two following types:
$$
\aligned
L(2p,\chi_f)&=\sum_{n=1}^{\infty}\frac{f(n)}{n^{2p}},\quad p\ge 1\ (r=2p),\text{ if }f(-1)=1.\\
L(2p-1,\chi_f)&=\sum_{n=1}^{\infty}\frac{f(n)}{n^{2p-1}},\quad p\ge 1\ (r=2p-1),\text{ if }f(-1)=-1.
\endaligned\tag1.1
$$

\flushpar An example of a
symmetric function mod $2m$ that will be used throughout this paper is the function $\chi_{2m}\: \{1,\dots ,2m-1\}\to \{0,\pm 1\}$, $m\ge 2$, such that
$$
\chi_{2m}(a)=\cases 1\text{ if } 1\le a\le m-1\\
-1\text{ if } m+1\le a\le 2m-1\\
0\text{ if } a=m
\endcases \tag1.2
$$
and $\chi_{2m}$ is extended by periodicity according to $(S_2)$. For $r$ odd, $\chi_{2m}(2m-a)=-\chi_{2m}(a)=(-1)^r\chi_{2m}(a)$ for all $a\in\{1,\dots ,2m-1\}$. 
Hence $\chi_{2m}$ is a symmetric function mod $2m$ for all $r$ odd. For $m=2$, $\chi_4(1)=1$, $\chi_4(3)=-1$, hence $\chi_4$ is a
homomorphism on $\bold Z_4^{\times}$, i.e., $\chi_4$ is a Dirichlet character mod 4. However, $\chi_{2m}$ is not a Dirichlet character for all $m\ge 3$. To see this,
(i) if $m=ab$, where $a>1$, $b>1$ then $\chi_{2m}(a)\chi_{2m}(b)=1\ne 0=\chi_{2m}(m)$; (ii) if $m\ge 3$ is prime then $m+1=2k$
and $\chi_{2m}(2)\chi_{2m}(k)=1\ne -1=\chi_{2m}(m+1)$ (in case (ii) both 2, $\frac{m+1}{2}\le m-1$ and $m+1< 2m-1$).\vskip.25cm

\flushpar The associated $L$-series for the symmetric function $\chi_{2m}$, $r$ odd, is
$$
\aligned
L(r,\chi_{2m})&=\sum_{n=1}^{\infty}\frac{\chi_{2m}(n)}{n^r}=\sum_{k=0}^{\infty}\sum_{a=1}^{2m-1}\frac{\chi_{2m}(a)}{(2km+a)^r}\\
&=\sum_{k=0}^{\infty}\left[\sum_{a=1}^{m-1}\frac{1}{(2km+a)^r}-\sum_{a=m+1}^{2m-1}\frac{1}{(2km+a)^r}\right]
\endaligned\tag1.3
$$
The $L$-series in Theorem 3.5, $r=3$, which include the examples above $L(3,\chi_8)$, $L(3,\chi_{12})$, are special cases of the above
general $L$-series (1.3) associated to the symmetric function $\chi_{2m}$ in (1.2). \vskip.25cm

\flushpar In \S 3 (3.5), we prove the following formula for the $L$-series (1.3) in case $r=3$.
$$
L(3,\chi_{2m})=\frac{\pi^3}{4m^3}\sum_{a=1}^{m-1}\frac{\sin(a\pi/m)}{(1-\cos(a\pi/m))^2}.\tag1.4
$$
Theorem 3.5 is proved by explicitly calculating the formula (1.4) for small $m\ge 2$, including in particular the calculations
above for $L(3,\chi_8)$, $L(3,\chi_{12})$. There exist similar trigonometric formulas for $L(2r+1,\chi_{2m})$ for all $r\ge1$
(cf \S5 (5.3) for $L(5,\chi_{2m}))$.\vskip.25cm

\flushpar The summand for $k=0$ in (1.3) for the $L$-series  $L(r,\chi_{2m})$, $r$ odd $\ge 3$, is the finite series
$$
1+\frac{1}{2^r}+\dots +\frac{1}{(m-1)^r}-\frac{1}{(m+1)^r}-\dots -\frac{1}{(2m-1)^r}.\tag1.5
$$
As $m\to \infty$ clearly the $+$ terms in (1.5) dominate in the sum $L(r,\chi_{2m})$ so that one obtains in the limit:
$$
\aligned
\zeta(2r+1)&=\sum_{n=1}^{\infty}\frac{1}{n^{2r+1}}\\
&=\lim_{m\to\infty}L(2r+1,\chi_{2m}) \text{ for all }r\ge 1.\endaligned\tag1.6
$$
The limit formula (1.6) shows the interest in symmetric functions for the study of $\zeta(2r+1)$, $r\ge 1$. Euler's famous calculation of $\zeta(2n)$ for all $n\ge 1$ led to the corresponding question of a
formula for $\zeta(2n+1)$, $n\ge 1$. However, to this day a closed-form formula for $\zeta(2n+1)$ remains unknown for any $n\ge 1$. Since Euler's time, formulas for the sums of some
series related to $\zeta(n)$ have been found in connection with: (i) the theory of trigonometric series (e.g.,\,Bromwich [1]), (ii) the theory of residues
in complex analysis
(e.g.,\,Sansone and Gerretsen [4]), and (iii) the theory of Dirichlet $L$-series associated to Dirichlet characters mod $N$ (Kato et al [2]). In this paper, (iii) above is generalized to the theory
of $L$-series associated to symmetric functions mod $N$. This generalization provides new formulas for sums of $L$-series associated to symmetric functions mod $N$, including $L$-series more closely related to $\zeta(3)$, $\zeta(5)$ (cf \S3 Theorem 3.5, \S4 Theorem 4.1, \S5 Theorem 5.2).

\flushpar{\bf \S2.1.} In order to state the main theoretical result Theorem 2.3 on $L$-series associated to symmetric functions we introduce and review the main properties of the auxiliary functions $h_r\: \bold C-\{1\}\to \bold C$ for all $r\ge 1$, employed by Kato et al [2, page 85] in connection with Dirichlet $L$-series.\vskip.25cm

Let $h_1(t)=\frac{1+t}{2(1-t)}$. For each integer $r\ge 2$, define $h_r(t)=\left(t\frac{d}{dt}\right)^rh_1(t)$. In particular,
$$
\aligned
h_2(t)&=t\frac{d}{dt}\left(\frac{1+t}{2(1-t)}\right)=\frac{t}{(1-t)^2}.\\
h_3(t)&=t\frac{d}{dt}(h_2(t))=\frac{t+t^2}{(1-t)^3}.\endaligned\tag2.1
$$
\flushpar{\bf Lemma 2.1} {\it Let $t=e^{2\pi ix}$, $x\in \bold C\setminus \bold Z$. Then $\cot \pi x=-2ih_1(t)=-2ih_1(e^{2\pi ix})$}.\vskip.25cm

\flushpar{\bf Proof.} $h_1(e^{2\pi ix})=\frac{e^{2\pi ix}+1}{2(1-e^{2\pi ix})}$. Hence $-2ih_1(e^{2\pi ix})=\frac{i(e^{2\pi ix}+1)}{e^{2\pi ix}-1}=\cot \pi x$.\qed\vskip.25cm

\flushpar Following Sansone and Gerretsen [4, page 145],
$$
\aligned
\pi\cot \pi x&=\frac{1}{x}+\sum_{n=1}^{\infty}\frac{2x}{x^2-n^2}=\frac12\sum_{n=-\infty}^{\infty}\frac{2x}{x^2-n^2}\\
&=\frac12\sum_{n=-\infty}^{\infty}\left(\frac{1}{x+n}+\frac{1}{x-n}\right),\quad x\in \bold C\setminus\bold Z.\endaligned\tag2.2
$$
Consequently, from Lemma 2.1
$$
h_1(t)=h_1(e^{2\pi ix})=-\frac12\cdot\frac{1}{2\pi i}\sum_{n=-\infty}^{\infty}\left(\frac{1}{x+n}+\frac{1}{x-n}\right),\quad x\in \bold C\setminus\bold Z.\tag2.3
$$
\flushpar{\bf Lemma 2.2} {\it For all} $r\ge 2$, $t=e^{2\pi ix}$,
$$
h_r(t)=h_r(e^{2\pi ix})=(r-1)!\left(\frac{-1}{2\pi i}\right)^r\sum_{n\in \bold Z}\frac{1}{(x+n)^r},\quad x\in \bold C\setminus\bold Z.
$$
\flushpar{\bf Proof.} Since $t=e^{2\pi ix}$ then $\dfrac{d}{dx}=\dfrac{dt}{dx}\dfrac{d}{dt}\Leftrightarrow      t\dfrac{d}{dt}=\dfrac{1}{2\pi i}\dfrac{d}{dx}$. Employing (2.3) one calculates
$$
\aligned
h_2(t)&=t\frac{d}{dt}(h_1(t))=\frac{1}{2\pi i}\frac{d}{dx}(h_1(t))\\
&=-\frac12\cdot\frac{1}{(2\pi i)^2}\sum_{n\in \bold Z}\left(-\frac{1}{(x+n)^2}-\frac{1}{(x-n)^2}\right)\\
\therefore\quad h_2(t)&=\frac{1}{(2\pi i)^2}\sum_{n\in\bold Z}\frac{1}{(x+n)^2}.\endaligned\tag2.4
$$
Thus (2.4) proves the lemma for $r=2$. By definition $h_{r+1}(t)=t\frac{d}{dt}(h_r(t))= \frac{1}{2\pi i}\frac{d}{dx}(h_r(t))$. Starting with (2.4) the
lemma for all $r\ge 2$ is proved by induction, where the sum of the series of term-by-term derivatives converges uniformly on compact subsets by the
Weierstrass $M$-test (comparison to the convergent series $\sum_{n\ge 1}\frac{1}{n^r}$, $r\ge 2$) \qed\vskip.25cm

\flushpar{\bf \S2.2.} The main result about $L$-series associated to symmetric functions is the following theorem which calculates $L$-series as a
finite sum in terms of the auxiliary functions $h_r(t)$. \vskip.25cm

\flushpar{\bf Theorem 2.3} {\it Let $r\ge 2$, $N\ge 2$. Let $\zeta_N=e^{\frac{2\pi i}{N}}$. Let $\chi(=\chi(N,r))\: \bold Z\to \bold C$ be a symmetric
function mod $N$, i.e., satisfying $(S_1)$, $(S_2)$ above. The associated $L$-series satisfies the following formula expressed in terms of the auxiliary function $h_r(t)$}:
$$
L(r,\chi)=\sum_{n=1}^{\infty}\frac{\chi(n)}{n^r}=\frac12\cdot\frac{1}{(r-1)!}\left(\frac{-2\pi i}{N}\right)^r\cdot\sum_{a=1}^{N-1}\chi(a)h_r(\zeta^a_N).
$$
\flushpar{\bf Proof of Theorem 2.3.} The right side of the formula for $L(r,\chi)$ in Theorem 2.3 is analyzed as follows. For all $r\ge 2$,
applying Lemma 2.2 to $t=\zeta_N^a =e^{2\pi ix}$, $x=a/N$, employing also the periodicity propert $(S_2)$ of the character $\chi$,
$$
\aligned
h_r(\zeta_N^a)&=(r-1)!\left(\frac{-1}{2\pi i}\right)^r\sum_{k\in\bold Z}\frac{1}{(\frac{a}{N}+k)^r},\quad a\in\{1,2,\dots ,N-1\}.\\
\therefore \sum_{a=1}^{N-1}\chi(a)h_r(\zeta_N^a)&=(r-1)!\left(\frac{-1}{2\pi i}\right)^r\sum_{k\in\bold Z}\sum_{a=1}^{N-1}\frac{N^r\chi(a)}{(kN+a)^r}\\
&=(r-1)!\left(\frac{-1}{2\pi i}\right)^r\sum_{k\in\bold Z}\sum_{a=1}^{N-1}\frac{N^r\chi(kN+a)}{(kN+a)^r}.
\endaligned\tag2.5
$$
We write $\sum_{k\in\bold Z}=\sum_{k\ge 0}+\sum_{k\le -1}$. From $(S_2)$, $\chi(kN)=0$ for all $k\in\bold Z$. Hence
$$
\sum_{k=0}^{\infty}\sum_{a=1}^{N-1}\frac{N^r\chi(kN+a)}{(kN+a)^r}=\sum_{k=0}^\infty\left[\sum_{kN<m<(k+1)N}\frac{N^r\chi(m)}{m^r}\right]=N^rL(r,\chi).\tag2.6
$$
Employing Remark 1.2, $\chi(-m)=(-1)^r\chi(m)$ for all $m\in \bold Z$. Again, since $\chi(kN)=0$ for all $k\in \bold Z$, one calculates the sum
$$
\aligned
\sum^{-\infty}_{k=-1}\sum_{a=1}^{N-1}\frac{N^r\chi(kN+a)}{(kN+a)^r}&=\sum^{\infty}_{k=1}\sum_{a=1}^{N-1}\frac{N^r\chi(-kN+a)}{(-kN+a)^r}=
\sum^{\infty}_{k=1}\sum_{a=1}^{N-1}\frac{N^r\chi(-(kN-a))}{(-(kN-a))^r}\\
&=\sum_{k=1}^{\infty}\sum_{a=1}^{N-1}
\frac{(-1)^r\cdot N^r\cdot\chi(kN-a)}{(-1)^r\cdot(kN-a)^r}\\
&=\sum_{k=0}^{\infty}\left[\sum_{kN<m<(k+1)N}\frac{N^r\chi(m)}{m^r}\right]=N^rL(r,\chi).
\endaligned\tag2.7
$$
Applying (2.6), (2.7) to (2.5) one obtains the formula,
$$
\aligned
\sum_{a=1}^{N-1}\chi(a)h_r(\zeta_N^a)&=(r-1)!\left(\frac{-1}{2\pi i}\right)^r\sum_{k\in\bold Z}\sum_{a=1}^{N-1}\frac{N^r\chi(kN+a)}{(kN+a)^r}\\
&=(r-1)!\left(\frac{-1}{2\pi i}\right)^r\left[\sum_{k=0}^{\infty}\sum_{a=1}^{N-1}\frac{N^r\chi(kN+a)}{(kN+a)^r}+
\sum_{k=-1}^{-\infty}\sum_{a=1}^{N-1}\frac{N^r\chi(kN+a)}{(kN+a)^r}\right]\\
&=(r-1)!\left(\frac{-1}{2\pi i}\right)^r\cdot 2N^rL(r,\chi).\endaligned\tag2.8
$$
Solving (2.8) for $L(r,\chi)$, one obtains
$$
L(r,\chi)=\frac12\cdot\frac{1}{(r-1)!}\left(\frac{-2\pi i}{N}\right)^r\cdot\sum_{a=1}^{N-1}\chi(a)h_r(\zeta_N^a)
$$
which completes the proof of Theorem 2.3. \qed
\vskip.25cm
\flushpar{\bf Historical Note 2.4.} In case $\chi\:\bold Z_N^{\times}\to \bold C^{\times}$ is a Dirichlet character mod $N$ then a formula analogous to
Theorem 2.3 for $L(r,\chi)$ is proved in Kato et al [2, Theorem 3.4, page 86], where these authors suppose also that
$\chi(-1)=(-1)^r$. As explained in Remark 1.3, Remark 1.4, \S1, a Dirichlet character $\chi$ is a special case of a symmetric function mod $N$. The contribution
of Theorem 2.3
above is that the symmetric function property of $\chi$ alone suffices to obtain the formula for $L(r,\chi)$ in Theorem 2.3, i.e., the Dirichlet character
hypothesis of Kato et al [2, p.\,86] on $\chi$
is not required. This theoretical point, which seems unrecognized in the $L$-series literature, was the initial inspiration for this
paper on $L$-series associated to symmetric functions.
The scope of $L$-series amenable to calculation by Theorem 2.3 is far greater than that of $L$-series based only on Dirichlet characters.\vskip.25cm

\flushpar{\bf \S 3.1. Analytic properties of the functions $h_r(t)$.} \vskip.25cm

\flushpar{\bf Lemma 3.1} (i) {\it If $r\ge 1$ is odd, then the image $h_r(S^1-\{1\})\subset \bold R i\subset \bold C$.}
(ii) {\it If $r\ge 2$ is even, then the image $h_r(S^1-\{1\})\subset \bold R$}. (iii) {\it For all $r\ge 1$, $t\in \bold S^1-\{1\}$,
$h_r(\overline t)(=h_r(\tfrac{1}{t}))=(-1)^rh_r(t)$. In particular if $r$ is odd then $h_r(-1)=0$.}\vskip.25cm

\flushpar{\bf Proof.} Let $t=e^{2\pi ix}\in S^1-\{1\}$, $x\in (0,1)\subset \bold R$. If $r\ge 2$, Lemmas 3.1(i)(ii) follow from the series expansion of $h_r(t)$ 
in Lemma 2.2, where $i^r\in\bold R$ if and
only if $r$ is even. Since $\overline t=\tfrac{1}{t}=e^{2\pi i(1-x)}$, $x\in (0,1)$, Lemma 3.1(iii) follows from Lemma 2.2 for $h_r(\overline t)$:
$$
\sum_{n\in \bold Z}\frac{1}{(1-x+n)^r}=\sum_{n\in \bold Z}\frac{1}{(-x+n)^r}=\sum_{n\in \bold Z}(-1)^r\frac{1}{(x-n)^r}=(-1)^r\sum_{n\in \bold Z}\frac{1}{(x+n)^r}.
$$
If $r=1$ the coefficient $i$ in the series expansion (2.3) of $h_1(t)$ proves Lemma 3.1(i) in this case. Also $h_1(\tfrac{1}{t})=\frac{1+\frac{1}{t}}{2(1-\frac{1}{t})}=-h_1(t)$,
which proves Lemma 3.1(iii) for $r=1$. \qed
\vskip.25cm
 \flushpar{\bf Corollary 3.2} {\it Let $\zeta_N=e^{\frac{2\pi i}{N}}$, $N\ge 2$. Then $h_r(\zeta_N^{N-a})=(-1)^rh_r(\zeta_N^a)$,
 for all $r\ge 1$, $a\in\{1,2,\dots ,N-1\}$.}\vskip.25cm

 \flushpar{\bf Proof.} Evidently $\zeta_N^a\in S^1-\{1\}$ and $\zeta_N^{N-a}=\overline{\zeta_N^a}=\frac{1}{\zeta_n^a}$ for all $a\in\{1,2,\dots ,N-1\}$. The corollary follows
 from Lemma 3.1(iii).  \qed\vskip.25cm

 \flushpar{\bf Corollary 3.3 (Refinement of Theorem 2.3)} {\it In Theorem 2.3 let $N=2q$, $q\ge 2$, $r=2p+1$, $p\ge 1$. Then}
 $$
L(r,\chi)=\frac{(-1)^{p+1}}{(r-1)!}\cdot\frac{\pi^r i}{q^r}\cdot\sum_{a=1}^{q-1}\chi(a)h_r(\zeta_{2q}^a).
$$
\flushpar{\bf Proof.} From property $(S_1)$ and Corollary 3.2, $\chi(2q-a)h_r(\zeta_{2q}^{2q-a})=\chi(a)h_r(\zeta_{2q}^a)$. Also from $(S_1)$, $\chi(q)=0$.
Consequently from Theorem 2.3,
$$
\aligned
L(r,\chi)&=\frac12\cdot\frac{1}{(r-1)!}\left(\frac{-\pi i}{q}\right)^r\cdot\left[\sum_{a=1}^{q-1}\chi(a)h_r(\zeta_{2q}^a)+\chi(2q-a)h_r(\zeta_{2q}^{2q-a})\right]\\
&=\frac12\cdot\frac{(-1)^{p+1}}{(r-1)!}\cdot\frac{\pi^ri}{q^r}\cdot2\sum_{a=1}^{q-1}\chi(a)h_r(\zeta_{2q}^a)\\
&=\frac{(-1)^{p+1}}{(r-1)!}\cdot\frac{\pi^ri}{q^r}\cdot\sum_{a=1}^{q-1}\chi(a)h_r(\zeta_{2q}^a).\qed
\endaligned
$$

\flushpar{\bf \S3.2. Trigonometric properties of the functions $h_r(t)$.} From Lemma 2.1 if $t=e^{2\pi ix}$, $x\in \bold C\setminus \bold Z$, then
$h_1(t)=\frac{-1}{2i}\cot (\pi x)=\frac{i}{2}\cot (\pi x)$. In addition from (2.4), $t\frac{d}{dt}=\frac{1}{2\pi i}\frac{d}{dx}$. Hence for all $r\ge 1$,
$$
h_{r+1}(t)=t\frac{d}{dt}(h_r(t))=\frac{1}{2\pi i}\frac{d}{dx}(h_r(t)),\quad t=e^{2\pi ix},\  x\in\bold C\setminus \bold Z.\tag3.1
$$
\flushpar{\bf Lemma 3.4} {\it If $t=e^{2\pi ix}$, $x\in \bold C\setminus \bold Z$},
\item{(i)} $h_2(t)=-\frac14\csc^2(\pi x)= -\frac14\frac{1}{\sin^2\pi x}=-\frac12\frac{1}{1-\cos 2\pi x}$.
\item{(ii)} $h_3(t)=-\frac{i}{2}\frac{\sin 2\pi x}{(1-\cos 2\pi x)^2}=-\frac{i}{4}\frac{\cos \pi x}{\sin^3\pi x}$. \vskip.25cm
\flushpar{\bf Proof.} Employing (3.1),
\item{(i)} $h_2(t)=\frac{1}{2\pi i}\frac{d}{dx}(h_1(t))=\frac{1}{2\pi i}\frac{d}{dx}(\frac{i}{2}\cot(\pi x))=-\frac14\csc^2\pi x=-\frac12\frac{1}{1-\cos 2\pi x}$.
\item{(ii)} $h_3(t)=\frac{1}{2\pi i}\frac{d}{dx}(h_2(t))=\frac{1}{2\pi i}\frac{d}{dx}(-\frac12\frac{1}{1-\cos 2\pi x})=
-\frac{i}{2}\frac{\sin 2\pi x}{(1-\cos 2\pi x)^2}=-\frac{i}{4}\frac{\cos \pi x}{\sin^3\pi x}$.\qed
\vskip.25cm

\flushpar{\bf \S3.3.} In this section we state and prove the main result Theorem 3.5 which calculates explicitly the series $L(3, \chi_{2m})$ for small values of $m$,
where $\chi_{2m}$, $m\ge 2$, is the symmetric function defined in (1.2). Applying Corollary 3.3 to $\chi_{2m}$, $N=2m$, $r=2p+1$, $p\ge 1$, the series $L(r,\chi_{2m})$ simplifies (recall from (1.2), $\chi_{2m}(a)=1$, $1\le a\le m-1$):
$$
\aligned
L(r,\chi_{2m})&=\sum_{n-1}^{\infty}\frac{\chi_{2m}(n)}{n^r}\\
&=\frac{(-1)^{p+1}}{(r-1)!}\cdot\frac{\pi^ri}{m^r}\cdot\sum_{a=1}^{m-1}\chi_{2m}(a)h_r(\zeta_{2m}^a)\\
&=\frac{(-1)^{p+1}}{(r-1)!}\cdot\frac{\pi^ri}{m^r}\cdot\sum_{a=1}^{m-1}h_r(\zeta_{2m}^a).
\endaligned\tag3.2
$$
In particular, in the case $r=3$ $(p=1)$,
$$
L(3,\chi_{2m})=\frac{\pi^3 i}{2m^3}\sum_{a=1}^{m-1}h_3(\zeta_{2m}^a).\tag3.3
$$
Employing Lemma 3.4(ii), one can replace the terms $h_3(\zeta_{2m}^a)$ by their corresponding trigonometric values, which turn out to be more convenient for calculations. Specifically, let $t=\zeta_{2m}^a=e^{\frac{2\pi ia}{2m}}=e^{2\pi ix}$, $x=\frac{a}{2m}$. Applying Lemma 3.4(ii),
$$
h_3(\zeta_{2m}^a)=-\frac{i}{2}\frac{\sin 2\pi x}{(1-\cos 2\pi x)^2}=-\frac{i}{2}\frac{\sin\frac{\pi a}{m}}{(1-\cos \frac{\pi a}{m})^2}.\tag3.4
$$
Employing (3.4), the formula (3.3) for $L(3,\chi_{2m})$ can be expressed in trigonometric terms,
$$
\aligned
L(3,\chi_{2m})&=\frac{\pi^3 i}{2m^3}\sum_{a=1}^{m-1}\left(-\frac{i}{2}\frac{\sin\frac{\pi a}{m}}{(1-\cos \frac{\pi a}{m})^2}\right)\\
&=\frac{\pi^3}{4m^3}\left[\frac{\sin \frac{\pi}{m}}{(1-\cos \frac{\pi}{m})^2}+\dots +\frac{\sin\tfrac{(m-1)\pi}{m}}{\left(1-\cos\tfrac{(m-1)\pi}{m}\right)^2}\right].
\endaligned\tag3.5
$$
We illustrate this formula in the simplest case $m=2$. According to \S1 (1.2), $\chi_4(1)=1$; $\chi_4(2)=0$; $\chi_4(3)=-1$, and it extends by periodicity to a symmetric function (same notation) $\chi_4\:\bold Z\to \{0,\pm 1\}$. Clearly $\chi_4\:\bold Z_4^{\times}\to \bold R^{\times}$ is a Dirichlet character (the unique Dirichlet character among the $\chi_{2m}$, $m\ge 2$). Its associated $L$-series is well-known:
$$
L(3,\chi_4)=1-\frac{1}{3^3}+\frac{1}{5^3}-\frac{1}{7^3}+\dots =\frac{\pi^3}{32}.
$$
Applying formula (3.5) to the case $m=2$ one obtains (there is only one summand in this case):
$$
L(3,\chi_4)=\frac{\pi^3}{4\cdot 8}\cdot \frac{\sin\pi/2}{(1-\cos\pi/2)^2}=\frac{\pi^3}{32}.
$$

\flushpar{\bf Theorem 3.5} {\it Employing (3.5) we calculate $L(3,\chi_{2m})$ for $m\in\{3,4, 6, 12\}$}.
$$
\aligned
L(3,\chi_6)&=\left(1+\frac{1}{2^3}-\frac{1}{4^3}-\frac{1}{5^3}\right)+\left(\frac{1}{7^3}+\frac{1}{8^3}-\frac{1}{10^3}-\frac{1}{11^3}\right)
 +\dots =\frac{5\pi^3\sqrt{3}}{243}.\\
L(3,\chi_8)&=\left(1+\frac{1}{2^3} +\frac{1}{3^3}-\frac{1}{5^3}-\frac{1}{6^3}-\frac{1}{7^3}\right)\\
&+\left(\frac{1}{9^3}+\frac{1}{10^3}+\frac{1}{11^3}-\frac{1}{13^3}-\frac{1}{14^3}-\frac{1}{15^3}\right)+\dots =\frac{\pi^3(6\sqrt{2}+1)}{256}.\\
L(3,\chi_{12})&=\left(1+\frac{1}{2^3}+\dots +\frac{1}{5^3}-\frac{1}{7^3}-\dots -\frac{1}{11^3}\right)\\
&+\left(\frac{1}{13^3}+\dots+\frac{1}{17^3}-\frac{1}{19^3}-\dots -\frac{1}{23^3}\right)
+\dots =\frac{\pi^3(20\sqrt{3}+261)}{7776}.\endaligned
$$
$$
\aligned
L(3,\chi_{24})&=\left(1+\frac{1}{2^3}+\dots +\frac{1}{11^3}-\frac{1}{13^3}-\dots -\frac{1}{23^3}\right)\\
&+\left(\frac{1}{25^3}+\dots +\frac{1}{35^3}-\frac{1}{37^3}-\dots -\frac{1}{47^3}\right)+\dots\\
&=\frac{\pi^3}{62,208}\bigg[(2484-828\sqrt{3})(2+\sqrt{3})^{1/2}+54\sqrt{2}+20\sqrt{3}+261\bigg].
\endaligned
$$

\flushpar{\bf Proof of Theorem 3.5.} The following ``double-angle'' formula is useful for the calculations, where in formula (3.6) we group together pairs of terms
involving $\theta$, $\pi-\theta$, such that $\theta=\frac{\pi a}{m}$, $1\le a\le m-1$:
$$
\frac{\sin\theta}{(1-\cos\theta)^2}+\frac{\sin(\pi-\theta)}{(1-\cos(\pi-\theta))^2}=\frac{4\sin\theta(3+\cos 2\theta)}{(1-\cos 2\theta)^2}, \quad \theta\in(0,\pi).\tag3.6
$$
\flushpar{(i)} Employing formula (3.5) for $m=3$, and also (3.7) for $\theta =\pi/3$,
$$
\aligned
L(3,\chi_6)&=\frac{\pi^3}{4\cdot 27}\left[\frac{\sin\pi/3}{(1-\cos\pi/3)^2}+\frac{\sin2\pi/3}{(1-\cos2\pi/3)^2}\right]\\
&=\frac{\pi^3}{4\cdot27}(\sin\pi/3)\left[\frac{4(3+\cos2\pi/3)}{(1-\cos2\pi/3)^2}\right]=\frac{\pi^3}{108}\cdot\frac{\sqrt{3}}{2}\cdot
\left[\frac{4(3-\frac12)}{(1+\frac12)^2}\right]=\frac{\pi^3 5\sqrt{3}}{243}.
\endaligned
$$
\flushpar{(ii)} Employing formula (3.5) for $m=4$, and also (3.7) for $\theta=\pi/4$,
$$
\aligned
L(3,\chi_8)&=\frac{\pi^3}{4\cdot 64}\left[\frac{\sin\pi/4}{(1-\cos\pi/4)^2}+\frac{\sin3\pi/4}{(1-\cos3\pi/4)^2}+\frac{\sin\pi/2}{(1-\cos\pi/2)^2}\right]\\
&=\frac{\pi^3}{4\cdot 64}\left[(\sin\pi/4)\left(\frac{4(3+\cos\pi/2)}{(1-\cos\pi/2)^2}\right)+1\right]=\frac{\pi^3}{256}(6\sqrt{2}+1).
\endaligned
$$
\flushpar(iii) Employing formula (3.5) for $m=6$, and also (3.7) for $\theta\in \{\pi/6, 2\pi/6=\pi/3\}$,
$$
\aligned
L(3,\chi_{12})&=\frac{\pi^3}{4\cdot 6^3}\left[\sum_{a=1}^5\frac{\sin(\pi a/6)}{(1-\cos(\pi a/6))^2}\right]\\
&=\frac{\pi^3}{4\cdot 6^3}\left[\sum_{a=1}^{a=2}(\sin\pi a/6)\left(\frac{(4(3+\cos 2\pi a/6)}{(1-\cos 2\pi a/6)^2}\right)+
\frac{\sin\pi/2}{(1-\cos\pi/2)^2}\right]\\
&=\frac{\pi^3}{4\cdot 6^3}\left[\frac12\left(\frac{4(3+\frac12)}{(1-\frac12)^2}\right)+\frac{\sqrt{3}}{2}\left(\frac{4(3-\frac12)}{(1+\frac12)^2}\right)+1\right]\\
&=\frac{\pi^3}{4\cdot 6^3}\left[28+\frac{20\sqrt{3}}{9}+1\right]=\frac{\pi^3}{7776}(20\sqrt{3}+261).
\endaligned
$$
\flushpar(iv) Employing formula (3.5) for $m=12$, and also (3.7) for $\theta\in \{\pi a/12\ \big |\  1\le a\le 5\}$
$$
\aligned
L(3,\chi_{24})&=\frac{\pi^3}{4\cdot 12^3}\left[\sum_{a=1}^{11}\frac{\sin(\pi a/12)}{(1-\cos(\pi a/12))^2}\right]\\
&=\frac{\pi^3}{4\cdot 12^3}\left[\sum_{a=1}^5(4\sin\pi a/12)\left(\frac{3+\cos2\pi a/12}{(1-\cos2\pi a/12)^2}\right)+\frac{\sin\pi/2}{(1-\cos\pi/2)^2}\right]\\
&=\frac{\pi^3}{6912}\left[\sum_{i=1}^5 A_a+1\right],\quad A_a=(\sin\pi a/12)\cdot\frac{4(3+\cos(2\pi a/12))}{(1-\cos(2\pi a/12))^2}, \quad1\le a\le 5.
\endaligned
$$
The calculations of $A_a$, $1\le a\le 5$ are as follows. Note that from (iii) above, $A_2=28$; $A_4=\frac{20\sqrt{3}}{9}$. From (ii) above, $A_3=6\sqrt{2}$.
$$
\aligned
A_1&=4(\sin\pi/12)\left[\frac{3+\cos\pi/6}{(1-\cos\pi/6)^2}\right]\\
&=4(\sin\pi/12)\left[\frac{3+\frac{\sqrt{3}}{2}}{\left(1-\frac{\sqrt{3}}{2}\right)^2}\right]=8(\sin\pi/12)\left[\frac{6+\sqrt{3}}{(2-\sqrt{3})^2}\right].\\
A_5&=4(\sin5\pi/12)\left[\frac{3+\cos5\pi/6}{(1-\cos5\pi/6)^2}\right]\\
&=4(\sin5\pi/12)\left[\frac{3-\frac{\sqrt{3}}{2}}{\left(1+\frac{\sqrt{3}}{2}\right)^2}\right]=8(\sin5\pi/12)\left[\frac{6-\sqrt{3}}{(2+\sqrt{3})^2}\right].
\endaligned
$$
Consequently from (iv) above,
$$
L(3,\chi_{24})=\frac{\pi^3}{6912}\left[A_1+A_5+28+6\sqrt{2}+\frac{20\sqrt{3}}{9}+1\right].\tag3.7
$$
We calculate the sum $A_1+A_5$. Recall $\cos2\theta=2\cos^2\theta-1=1-2\sin^2\theta$. For $\theta=\pi/12$,
$$
\sin\pi/12=\frac{(2-\sqrt{3})^{1/2}}{2}\ ;\quad \sin5\pi/12=\cos\pi/12=\frac{(2+\sqrt{3})^{1/2}}{2}.\tag3.8
$$
Employing (3.8) and $A_1,A_5$ above, noting also that $(2+\sqrt{3})(2-\sqrt{3})=1$,
$$
\aligned
A_1+A_5&= 8\frac{(2-\sqrt{3})^{1/2}}{2}\left[\frac{6+\sqrt{3}}{(2-\sqrt{3})^2}\right]+8\frac{(2+\sqrt{3})^{1/2}}{2}\left[\frac{6-\sqrt{3}}{(2+\sqrt{3})^2}\right]\\
&=\frac{4(6+\sqrt{3})}{(2-\sqrt{3})^{3/2}}+\frac{4(6-\sqrt{3})}{(2+\sqrt{3})^{3/2}}\\
&=4(6+\sqrt{3})(2+\sqrt{3})^{3/2}+4(6-\sqrt{3})(2-\sqrt{3})^{3/2}\\
&=4(15+8\sqrt{3})(2+\sqrt{3})^{1/2}+4(15-8\sqrt{3})(2-\sqrt{3})^{1/2}\\
&=4(2+\sqrt{3})^{1/2}\left[15+8\sqrt{3}+(15-8\sqrt{3})(2-\sqrt{3})\right]\\
&=(276-92\sqrt{3})(2+\sqrt{3})^{1/2}.
\endaligned\tag3.9
$$
Employing (3.7), (3.9) one calculates
$$
\aligned
L(3,\chi_{24})&=\frac{\pi^3}{6912}\left[(276-92\sqrt{3})(2+\sqrt{3})^{1/2}+29+6\sqrt{2}+\frac{20\sqrt{3}}{9}\right]\\
&=\frac{\pi^3}{62,208}\left[(2484-828\sqrt{3})(2+\sqrt{3})^{1/2}+54\sqrt{2}+20\sqrt{3}+261\right].\qed
\endaligned
$$

\flushpar{\bf \S4.1.} Let $\zeta^{\text{odd}}(r)=\sum_{k\ge 0}\frac{1}{(2k+1)^r}$, $r\ge 2$. We modify $\chi_{2m}$ to obtain a new symmetric
function whose values are zero on the even integers, and whose associated $L$-series is adapted to $\zeta^{\text{odd}}(r)$, $r$ odd $\ge 3$.
Let $f_{4m}\:\{1,2\dots ,4m-1\}\to \{0,\pm 1\}$, $m\ge 1$, such that
$$
f_{4m}(a)=\cases 0\text{ if }a\text{ is even}\\
1\text{ if $a$ is odd},\quad 1\le a\le 2m-1\\
-1\text{ if $a$ is odd},\quad 2m+1\le a\le 4m-1.
\endcases\tag4.1
$$
For $r$ odd, $f_{4m}(4m-a)=-f_{4m}(a)=(-1)^rf_{4m}(a)$, $1\le a\le 4m-1$. Hence for all $r$ odd, $m\ge 1$, $f_{4m}$ is a symmetric function mod $4m$ which extends by periodicity to a function (same notation) $f_{4m}\: \bold Z\to \{0,\pm 1\}$ such that $f_{4m}(a)=0$ if $a$ is even; $f_{4m}(4km+a)=f_{4m}(a)$ for all $k\in \bold Z$, $a\in\{1,2,\dots ,4m-1\}$; $f_{4m}(-k)=(-1)^rf_{4m}(k)$ for all $k\in\bold Z$.\vskip.25cm

\flushpar For all $m\ge 1$, $r$ odd, the $L$-series associated to the symmetric function $f_{4m}$ is
$$
\aligned
L(r,f_{4m})&=\sum_{n=1}^{\infty}\frac{f_{4m}(n)}{n^r}=\sum_{k=0}^{\infty}\sum_{a=1}^{4m-1}\frac{f_{4m}(a)}{(4km+a)^r}\\
&=\sum_{k=0}^{\infty}\left[\sum\Sb 1\le a\le 2m-1\\ a\text{ odd}\endSb\frac{1}{(4km+a)^r}-\sum\Sb 2m+1\le a\le 4m-1\\ a\text{ odd}\endSb\frac{1}{(4km+a)^r}\right]
\endaligned\tag4.2
$$
From (4.2), one obtains, in a similar way to formula \S1, (1.6),
$$
\aligned
\zeta^{\text{odd}}(2r+1)&=\sum_{n=0}^{\infty}\frac{1}{(2n+1)^{2r+1}}\\
&=\lim_{m\to \infty}L(2r+1,f_{4m}),\text{ for all }r\ge 1.
\endaligned\tag4.3
$$
\flushpar{\bf \S4.2. Calculation of $L(2r+1,f_{4m})$.} Applying Corollary 3.3 to $f_{4m}$, $m\ge 1$, in the case $N=4m$, $q=2m$, $r=2p+1$, $p\ge 1$, one obtains
$$
L(r,f_{4m})=\frac{(-1)^{p+1}}{(r-1)!}\cdot \frac{\pi^r i}{(2m)^r}\cdot\sum\Sb 1\le a\le 2m-1\\ a\text{ odd}\endSb  h_r(\zeta_{4m}^a).\tag4.4
$$
In particular, in the case $r=3$ $(p=1)$,
$$
\aligned
L(3,f_{4m}&)=\frac{\pi^3i}{2(2m)^3}\cdot\sum\Sb 1\le a\le 2m-1\\ a\text{ odd}\endSb  h_3(\zeta_{4m}^a)
=\frac{\pi^3i}{16m^3}\cdot\sum\Sb 1\le a\le 2m-1\\ a\text{ odd}\endSb  h_3(\zeta_{4m}^a).
\endaligned\tag4.5
$$
Let $t=\zeta_{4m}^a=e^{\frac{2\pi ia}{4m}}=e^{2\pi ix}$, $x=\frac{a}{4m}$. Applying Lemma 3.4(ii),
$$
h_3(\zeta_{4m}^a)=-\frac{i}{2}\frac{\sin 2\pi x}{(1-\cos 2\pi x)^2}=-\frac{i}{2}\frac{\sin\frac{\pi a}{2m}}{(1-\cos \frac{\pi a}{2m})^2}.\tag4.6
$$
Employing (4.6), the formula (4.5) for $L(3,f_{4m})$ can be expressed in trigonometric terms,
$$
\aligned
L(3,f_{4m})&=\frac{\pi^3 i}{16m^3}\cdot\sum\Sb 1\le a\le 2m-1\\ a\text{ odd}\endSb \left(-\frac{i}{2}\frac{\sin\frac{\pi a}{2m}}{(1-\cos \frac{\pi a}{2m})^2}\right)\\
&=\frac{\pi^3}{32m^3}\left[\frac{\sin \frac{\pi}{2m}}{(1-\cos \frac{\pi}{2m})^2}+\frac{\sin \frac{3\pi}{2m}}{(1-\cos \frac{3\pi}{2m})^2}+\dots +\frac{\sin\tfrac{(2m-1)\pi}{2m}}{\left(1-\cos\tfrac{(2m-1)\pi}{2m}\right)^2}\right].
\endaligned\tag4.7
$$
We illustrate this formula in the simplest case $m=1$. According to (4.1), $f_4(1)=1$, $f_4(2)=0$, $f_4(3)=-1$. Hence $f_4=\chi_4\:\bold Z\to \{0,\pm 1\}$. Consequently, applying formula (4.7) to the case $m=1$ (note that there is only one summand in this case), one confirms the calculation of $L(3,\chi_4)$ made prior to the statement of \S3, Theorem 3.5:
$$
\aligned
L(3,f_4)&=1-\frac{1}{3^3}+\frac{1}{5^3}-\frac{1}{7^3}+\dots\\
&=\frac{\pi^3}{32}\cdot \frac{\sin\pi/2}{(1-\cos\pi/2)^2}=\frac{\pi^3}{32}.\endaligned\tag4.8
$$
{\bf Theorem 4.1} {\it Employing (4.7) we calculate $L(3,f_{4m})$ for $m\in \{2,3,6\}$}.
$$
\aligned
L(3,f_8)&=\left(1+\frac{1}{3^3}-\frac{1}{5^3}-\frac{1}{7^3}\right)+\left(\frac{1}{9^3}+\frac{1}{11^3}-\frac{1}{13^3}-\frac{1}{15^3}\right)
 +\dots =\frac{3\pi^3\sqrt{2}}{128}.\\
L(3,f_{12})&=\left(1+\frac{1}{3^3} +\frac{1}{5^3}-\frac{1}{7^3}-\frac{1}{9^3}-\frac{1}{11^3}\right)\\
&+\left(\frac{1}{13^3}+\frac{1}{15^3}+\frac{1}{17^3}-\frac{1}{19^3}-\frac{1}{21^3}-\frac{1}{23^3}\right)+\dots =\frac{29\pi^3}{864}.\\
L(3,f_{24})&=\left(1+\frac{1}{3^3}+\frac{1}{5^3}+\dots +\frac{1}{11^3}-\frac{1}{13^3}-\frac{1}{15^3}-\dots -\frac{1}{23^3}\right)\\
&+\left(\frac{1}{25^3}+\frac{1}{27^3}+\dots +\frac{1}{35^3}-\frac{1}{37^3}-\frac{1}{39^3}-\dots -\frac{1}{47^3}\right)+\dots\\
&=\frac{\pi^3}{62,208}\bigg[(2484-828\sqrt{3})(2+\sqrt{3})^{1/2}+54\sqrt{2}\bigg].
\endaligned
$$
\flushpar{\bf Proof of Theorem 4.1.} The following lemma relates $L(r,f_{4m})$ to the $L$-series $L(r,\chi_{4m})$, $L(r,\chi_{2m})$. The calculations of the $L$-series in Theorem 3.5 are then employed to calculate $L(3,f_{4m})$  for small $m$.\vskip.25cm
\flushpar{\bf Lemma 4.2} {\it For all $r\text{ odd}\ge 3$, $m\ge 2$,} $L(r,f_{4m})=L(r,\chi_{4m})-\frac{1}{2^r}L(r,\chi_{2m})$.\vskip.25cm
\flushpar{\bf Proof of Lemma 4.2.} Employing formulas (1.3), (4.2) one calculates,
$$
\aligned
L(r,\chi_{4m})&=\sum_{n\ge 1}\frac{\chi_{4m}(n)}{n^r}\\
&=\sum_{n=1}^{2m-1}\frac{1}{n^r}-\sum_{n=2m+1}^{4m-1}\frac{1}{n^r}+\text{ etc.}\\
\endaligned
$$
$$\aligned
\therefore\ L(r,\chi_{4m})&=\sum\Sb 1\le n\le 2m-1\\n\text{ odd}\endSb\frac{1}{n^r}-\sum\Sb 2m+1\le n\le 4m-1\\n\text{ odd}\endSb\frac{1}{n^r}+\text{ etc.}\\
&+\sum_{n=1}^{m-1}\frac{1}{(2n)^r}-\sum_{n=m+1}^{2m-1}\frac{1}{(2n)^r}+\text{ etc.}\\
\endaligned
$$
Hence $L(r,\chi_{4m})=L(r,f_{4m})+\frac{1}{2^r}L(r,\chi_{2m})\Leftrightarrow\ L(r,f_{4m})=L(r,\chi_{4m})-\frac{1}{2^r}L(r,\chi_{2m})$.\qed\vskip.25cm

\flushpar Returning to the proof of Theorem 4.1, we use Lemma 4.2 in the case $r=3$.\vskip.25cm
\flushpar (i) Applying Theorem 3.5, formula (4.8) and Lemma 4.2 for $m=2$,
$$
\aligned
L(3,f_{8})&=L(3,\chi_8)-\frac18\cdot L(3,\chi_4)\\
&=\frac{\pi^3(6\sqrt{2}+1)}{256}-\frac18\cdot\frac{\pi^3}{32}=\frac{3\pi^3\sqrt{2}}{128}.
\endaligned
$$
Alternatively, applying the trigonometric formula (4.7) for $m=2$ and (3.7) for $\theta =\pi/4$,
$$
\aligned
L(3,f_8)&=\frac{\pi^3}{32\cdot 8}\left[\frac{\sin\pi/4}{(1-\cos\pi/4)^2}+\frac{\sin 3\pi/4}{(1-\cos3\pi/4)^2}\right]\\
&=\frac{\pi^3}{256}\left[(\sin\pi/4)\frac{4(3+\cos\pi/2)}{(1-\cos\pi/2)^2}\right]=\frac{\pi^3}{256}\cdot\frac{12}{\sqrt{2}}=\frac{3\pi^3\sqrt{2}}{128}.\endaligned
$$
\flushpar (ii) Applying Theorem 3.5 and Lemma 4.2 for $m=3$,
$$
\aligned
L(3,f_{12})&=L(3,\chi_{12})-\frac18 \cdot L(3,\chi_6)\\
&=\frac{\pi^3(20\sqrt{3}+261)}{7776}-\frac18\cdot \frac{5\pi^3\sqrt{3}}{243}\ \left(=-\frac{20\pi^3\sqrt{3}}{7776}\right)\\
&=\frac{261\pi^3}{7776}=\frac{29\pi^3}{864}.
\endaligned
$$
(iii) Applying Theorem 3.5 and Lemma 4.2 for $m=6$,
$$
\aligned
L(3,f_{24})&=L(3,\chi_{24})-\frac18 \cdot L(3,\chi_{12})\\
&=\frac{\pi^3}{62208}\left[(2484-828\sqrt{3})(2+\sqrt{3})^{1/2}+54\sqrt{2}+20\sqrt{3}+261\right]\\
&-\frac18\cdot\frac{\pi^3(20\sqrt{3}+261)}{7776}\ \left(=-\frac{\pi^3(20\sqrt{3}+261)}{62,208}\right).\\
\endaligned
$$
Hence $L(3,f_{24})=\frac{\pi^3}{62,208}\left[(2484-828\sqrt{3})(2+\sqrt{3})^{1/2}+54\sqrt{2}\right]$.\qed\vskip.25cm
\flushpar{\bf Remark 4.3.} The formula for $L(3,f_8)$ in Theorem 4.1 was calculated also in Bromwich [1, p. 364] by other means using trigonometric series.\vskip.25cm

\flushpar{\bf Remark 4.4.} The symmetric function $f_{4m}$ mod $4m$ is a Dirichlet character if and only if $m\in\{1,2\}$. Indeed, from (4.1), $f_4(1)=1$, $f_4(3)=-1$. Hence $f_4$ is a homomorphism on $\bold Z_4^{\times}$. Similarly employing (4.1), $f_8(1)=f_8(3)=1$; $f_8(5)=f_8(7)=-1$. Hence $f_8$ is a homomorphism on $\bold Z_8^{\times}$. We show however that $f_{4m}$ is not a homomorphism on $\bold Z_{4m}^{\times}$ for all $m\ge 3$. Referring to (4.1) we consider the three cases $2m+1\equiv a$ mod 3.\vskip.25cm

\flushpar (i) If $2m+1=3p$ then $f_{4m}(3)\cdot f_{4m}(p)=1\ne -1=f_{4m}(2m+1)$; (ii) If $2m+1=3q+1$ then $2m+3=3(q+1)$,
hence $f_{4m}(3)\cdot f_{4m}(q+1)=1\ne -1=f_{4m}(2m+3)$; (iii) If $2m+1=3s+2$ then $2m+5=3(s+2)$, hence $f_{4m}(3)\cdot f_{4m}(s+2)=1\ne -1=f_{4m}(2m+5)$.\vskip.25cm

\flushpar To justify these calculations, note that if $m\ge 3$ then: in (i) both $3, \frac{2m+1}{3}\le 2m-1$ and $2m+1\le 4m-1$; in (ii) both $3,
\frac{2m+3}{3}\le 2m-1$ and $2m+3\le 4m-1$; in (iii) both $3,\frac{2m+5}{3}\le 2m-1$ and $2m+5\le 4m-1$. Consequently, in all cases (i), (ii), (iii)
the symmetric function $f_{4m}$ is not a Dirichlet character for all $m\ge 3$.\vskip.25cm

\flushpar{\bf \S5.1. $L$-series associated to $\zeta(5)$.} In this section we state and prove Theorem 5.2 which calculates the series $L(5,\chi_{2m})$ for
small values of $m$, analogous to the calculations of $L(3,\chi_{2m})$ in Theorem 3.5, \S3. According to formula (3.2) for $r=5$ $(p=2)$,
$$
L(5,\chi_{2m})=\frac{(-1)^3}{4!}\frac{\pi^5i}{m^5}\sum_{a=1}^{m-1}h_5(\zeta_{2m}^a).\tag5.1
$$
The following lemma calculates the auxiliary function $h_5(t)$ required for formula (5.1).\vskip.25cm
\flushpar{\bf Lemma 5.1} Let $t=e^{2\pi ix}$, $x\in \bold C\setminus\bold Z$.
\item{(i)} $h_4(t)=\frac{1+2\cos^2\pi x}{8\sin^4\pi x}$.
\item{(ii)} $h_5(t)=\frac{i}{4}\frac{2\cos\pi x+\cos^3\pi x}{\sin^5\pi x}=\frac{i}{2}\frac{\cot\pi x(5+\cos2\pi x)}{(1-\cos2\pi x)^2}.$\vskip.25cm

\flushpar{\bf Proof.} From Lemma 3.4, $h_3(t)=-\frac{i}{4}\frac{\cos\pi x}{\sin^3\pi x}$. Recall (3.1): $t\frac{d}{dt}=\frac{1}{2\pi i}\frac{d}{dx}$. One calculates,
$$
h_4(t)=t\frac{d}{dt}(h_3(t))=\frac{1}{2\pi i}\cdot-\frac{i}{4}\cdot\frac{d}{dx}\left[\frac{\cos\pi x}{\sin^3\pi x}\right]=\frac{1+2\cos^2\pi x}{8\sin^4\pi x}.
$$
$$
\aligned
h_5(t)&=\frac{1}{2\pi i}\frac{d}{dx}(h_4(t))=\frac{1}{2\pi i}\frac{d}{dx}\left[\frac{1+2\cos^2\pi x}{8\sin^4\pi x}\right]\\
&=\frac{i}{4}\cdot\frac{2\cos\pi x+\cos^3\pi x}{\sin^5\pi x}=\frac{i}{2}\cdot\cot\pi x\cdot\frac{5+\cos2\pi x}{(1-\cos2\pi x)^2}.
\endaligned
$$
Employing Lemma 5.1, one can replace the terms $h_5(\zeta_{2m}^a)$ in (5.1) by their corresponding trigonometric values. Let
$t=\zeta_{2m}^a=e^{\frac{2\pi ia}{2m}}=e^{2\pi ix}$, $x=\frac{\pi a}{2m}$. Applying Lemma 5.1(ii),

$$
h_5(\zeta_5^a)=\frac{i}{2}\cdot\cot\frac{\pi a}{2m}\cdot\frac{5+\cos\pi a/m}{(1-\cos\pi a/m)^2}.\tag5.2
$$
Employing (5.2) the formula (5.1) for $L(5,\chi_{2m})$ can be expressed in trigonometric terms,
$$
\aligned
L(5,\chi_{2m})&= \frac{-\pi^5i}{4!\cdot m^5}\cdot\frac{i}{2}\sum_{a=1}^{m-1}\cot(\pi a/2m)\cdot\frac{5+\cos\pi a/m}{(1-\cos\pi a/m)^2}\\
&=\frac{\pi^5}{2\cdot4!\cdot m^5}\sum_{a=1}^{m-1}\cot(\pi a/2m)\cdot\frac{5+\cos\pi a/m}{(1-\cos\pi a/m)^2}.
\endaligned\tag5.3
$$
We illustrate formula (5.3) in the simplest case $m=2$ (there is only one summand in this case).
According to \S1 (1.2), $\chi_4$ is the Dirichlet character $\chi_4(1)=1$, $\chi_4(2)=0$, $\chi_4(3)=-1$ which extends by periodicity to a symmetric function
mod 4, $\chi_4\: \bold Z\to \{0,\pm 1\}$, whose associated $L$-series is the alternating series
$$
L(5,\chi_4)=1-\frac{1}{3^5}+\frac{1}{5^5}-\frac{1}{7^5}+\cdots\ .
$$
Applying formula (5.3) to the case $m=2$ one obtains
$$
\aligned
\sum_{n=0}^{\infty}\frac{(-1)^n}{(2n+1)^5}&=\frac{\pi^5}{2\cdot 4!\cdot 2^5}\cdot\cot(\pi/4)\cdot\frac{5+\cos\pi/2}{(1-\cos\pi/2)^2}\\
&=\frac{5\pi^5}{2\cdot4!\cdot 2^5}=\frac{5\pi^5}{1536}.
\endaligned
$$
This result is a special case ($p=2$; $E_4=5$) of the classical formula,
$$
\sum_{n=0}^{\infty}\frac{(-1)^n}{(2n+1)^{2p+1}}=\frac12\cdot\left(\frac{\pi}{2}\right)^{2p+1}\frac{E_{2p}}{(2p)!},
$$
proved in Sansone and Gerretsen [4, pp. 86, 144] using the theory of residues in complex analysis, where the coefficients $E_{2p}$ are
known as Euler coefficients. \vskip.25cm

\flushpar{\bf Theorem 5.2} Employing formula (5.3) we calculate $L(5,\chi_{2m})$, $m\in\{6,8,12,24\}$.
$$
\aligned
L(5,\chi_6)&=\left(1+\frac{1}{2^5}-\frac{1}{4^5}-\frac{1}{5^5}\right)+\left(\frac{1}{7^5}+\frac{1}{8^5}-\frac{1}{10^5}-\frac{1}{11^5}\right)
 +\dots =\frac{\pi^5 17\sqrt{3}}{8748}.\\
L(5,\chi_8)&=\left(1+\frac{1}{2^5} +\frac{1}{3^5}-\frac{1}{5^5}-\frac{1}{6^5}-\frac{1}{7^5}\right)\\
&+\left(\frac{1}{9^5}+\frac{1}{10^5}+\frac{1}{11^5}-\frac{1}{13^5}-\frac{1}{14^5}-\frac{1}{15^5}\right)+\dots =\frac{\pi^5(5+114\sqrt{2})}{49,152}.\\
L(5,\chi_{12})&=\left(1+\frac{1}{2^5}+\dots +\frac{1}{5^5}-\frac{1}{7^5}-\dots -\frac{1}{11^5}\right)\\
&+\left(\frac{1}{13^5}+\dots+\frac{1}{17^5}-\frac{1}{19^5}-\dots -\frac{1}{23^5}\right)
+\dots =\frac{\pi^5(3675+68\sqrt{3})}{1,119,744}.\endaligned
$$
$$
\aligned
L(5,\chi_{24})&=\left(1+\frac{1}{2^5}+\dots +\frac{1}{11^5}-\frac{1}{13^5}-\dots -\frac{1}{23^5}\right)\\
&+\left(\frac{1}{25^5}+\dots +\frac{1}{35^5}-\frac{1}{37^5}-\dots -\frac{1}{47^5}\right)+\dots\\
&=\frac{\pi^5}{35,831,808}\bigg[(143,460-47,820\sqrt{3})(2+\sqrt{3})^{1/2}+342\sqrt{2}+68\sqrt{3}+3675\bigg].
\endaligned
$$
\flushpar{\bf Proof.} Analogous to the proof of Theorem 3.5, the following double-angle formula is useful for the calculations, where in formula (5.3) we we group
together pairs of terms involving $\alpha$, $\frac{\pi}{2}-\alpha$, such that $\alpha=\frac{\pi a}{2m}$, $1\le a\le m-1$.\vskip.25cm

\flushpar{\bf Lemma 5.3} Let $\theta\in(0,\pi/2)$. Then
$$
\aligned
&\frac{\cot\theta[5+\cos2\theta]}{(1-\cos2\theta)^2}+\frac{\cot(\frac{\pi}{2}-\theta)[5+\cos2(\frac{\pi}{2}-\theta)]}{(1-\cos2(\frac{\pi}{2}-\theta))^2}\\
&=\frac{2(5+18\cos^22\theta+\cos^42\theta)}{\sin^52\theta}=\frac{2}{\sin2\theta}\cdot\frac{57+38\cos4\theta+\cos^24\theta}{(1-\cos4\theta)^2}.
\endaligned
$$
\flushpar{\bf Proof of Lemma 5.3.} The left side of the first equality is
$$
\aligned
&\frac{\cot\theta[5+\cos2\theta]}{(1-\cos2\theta)^2}+\frac{\tan\theta[5-\cos2\theta]}{(1+\cos2\theta)^2}\\
&=\frac{1}{(1-\cos^22\theta)^2}\left[\cot\theta(5+\cos2\theta)(1+\cos2\theta)^2+\tan\theta(5-\cos2\theta)(1-\cos2\theta)^2\right]\\
&=\frac{1}{\sin^42\theta}\left[(5+7\cos^22\theta)(\cot\theta+\tan\theta)+(11\cos2\theta+\cos^32\theta)(\cot\theta-\tan\theta)\right]
\endaligned
$$
Note that $\cot\theta+\tan\theta=\frac{2}{\sin2\theta}$; $\cot\theta-\tan\theta=\frac{2\cos2\theta}{\sin2\theta}$. Simplifying, the first
equality of the Lemma is proved. The second equality is proved by the usual double-angle formulas.\qed
\vskip.25cm
\flushpar\item{(i)} Employing formula (5.3) for $m=3$, and also Lemma 5.3 for $\theta =\pi/6$,
$$
\aligned
L(5,\chi_6)&=\frac{\pi^5}{2\cdot4!\cdot3^5}\left[\frac{\cot(\pi/6)[5+\cos\pi/3]}{(1-\cos\pi/3)^2}+\frac{\cot(\pi/3)[5+\cos2\pi/3]}{(1-\cos2\pi/3)^2}\right]\\
&=\frac{\pi^5}{2\cdot4!\cdot3^5}\left[\frac{2(5+18\cos^2\pi/3+\cos^4\pi/3)}{\sin^5\pi/3}\right]\\
&=\frac{\pi^5}{48\cdot3^5}\left[\frac{2(5+\frac{18}{4}+\frac{1}{16})}{\frac{9\sqrt{3}}{32}}\right]=\frac{\pi^517\sqrt{3}}{8748}.
\endaligned
$$
\flushpar\item{(ii)} Employing formula (5.3) for $m=4$, and also Lemma 5.3 for $\theta =\pi/8$,
$$
\aligned
L(5,\chi_8)&=\frac{\pi^5}{2\cdot4!\cdot 4^5}\sum_{a=1}^3\frac{\cot(\pi a/8)[5+\cos\pi a/4]}{(1-\cos\pi a/4)^2}\\
&=\frac{\pi^5}{2\cdot4!\cdot 4^5}\left[\frac{2(5+18\cos^2\pi/4+\cos^4\pi/4)}{\sin^5\pi/4}+\frac{\cot(\pi/4)(5+\cos\pi/2)}{(1-\cos\pi/2)^2}\right]\\
&=\frac{\pi^5}{2\cdot4!\cdot 4^5}\left[\frac{2(5+18/2+1/4)}{\frac{1}{4\sqrt{2}}}+5\right]=\frac{\pi^5(114\sqrt{2}+5)}{49,152}.
\endaligned
$$
\flushpar\item{(iii)} Employing formula (5.3) for $m=6$, and also Lemma 5.3 for $\theta\in\{\pi/6,2\pi/6=\pi/3\}$,
$$
\aligned
L(5,\chi_{12})&=\frac{\pi^5}{2\cdot4!\cdot 6^5}\sum_{a=1}^5\frac{\cot(\pi a/12)[5+\cos\pi a/6]}{(1-\cos\pi a/6)^2}\\
&=\frac{\pi^5}{2\cdot4!\cdot 6^5}\left[\sum_{a=1}^2\frac{2(5+18\cos^2\pi a/6+\cos^4\pi a/6)}{\sin^5\pi a/6}+\frac{\cot(\pi/4)(5+\cos\pi/2)}{(1-\cos\pi/2)^2}\right]\\
&=\frac{\pi^5}{2\cdot4!\cdot 6^5}\left[\frac{2(5+18(3/4)+9/16)}{\frac{1}{32}}+\frac{2(5+18(1/4)+1/16)}{\frac{9\sqrt{3}}{32}}+5\right]\\
&=\frac{\pi^5}{2\cdot4!\cdot 6^5}\left[1220+\frac{68\sqrt{3}}{3}+5\right]=\frac{3675+68\sqrt{3}}{1,119,744}.
\endaligned
$$
\flushpar\item{(iv)} Employing formula (5.3) for $m=12$, and also Lemma 5.3 for $\theta\in\{\pi a/24
\mid 1\le a\le 5\}$,
$$
\aligned
L(5,\chi_{24})&=\frac{\pi^5}{2\cdot4!\cdot 12^5}\sum_{a=1}^{11}\frac{\cot(\pi a/24)[5+\cos\pi a/12]}{(1-\cos\pi a/12)^2}\\
&=\frac{\pi^5}{2\cdot4!\cdot 12^5}\left[\sum_{a=1}^5\frac{2(5+18\cos^2\frac{\pi a}{12}+\cos^4\frac{\pi a}{12})}{\sin^5\frac{\pi a}{12}}+
\frac{\cot(\pi/4)(5+\cos\pi/2)}{(1-\cos\pi/2)^2}\right]\\
&=\frac{\pi^5}{2\cdot4!\cdot 12^5}\left[\sum_{a=1}^5B_a+5\right],\quad
B_a=\frac{2(5+18\cos^2\frac{\pi a}{12}+\cos^4\frac{\pi a}{12})}{\sin^5\frac{\pi a}{12}}, \quad 1\le a\le 5.
\endaligned
$$
Note that from (iii) above, $B_2=1220$; $B_4=\frac{68\sqrt{3}}{3}$. From (ii) above, $B_3=114\sqrt{2}$. Consequently
$$
L(5,\chi_{24})=\frac{\pi^5}{2\cdot4!\cdot 12^5}\cdot\left[B_1+B_5+1220+\frac{68\sqrt{3}}{3} +114\sqrt{2}+5\right].\tag5.4
$$
Applying Lemma 5.3 (second equality), we calculate $B_1+B_5$:
$$
\aligned
B_1+B_5&=\frac{2(57+38\cos\pi/6+\cos^2\pi/6)}{\sin(\pi/12)(1-\cos\pi/6)^2}+\frac{2(57+38\cos5\pi/6+\cos^25\pi/6)}{\sin(5\pi/12)(1-\cos5\pi/6)^2}\\
&=\frac{2}{\sin\pi/12}\cdot\left[\frac{57+38\frac{\sqrt{3}}{2}+\frac34}{(1-\sqrt{3}/2)^2}\right]+
\frac{2}{\cos\pi/12}\cdot\left[\frac{57-38\frac{\sqrt{3}}{2}+\frac34}{(1+\sqrt{3}/2)^2}\right].
\endaligned
$$
Employing (3.9) for the terms $\sin\pi/12$, $\cos\pi/12$, noting also that $(2+\sqrt{3})(2-\sqrt{3})=1$, one calculates,
$$
\aligned
B_1+B_5&=\frac{4}{(2-\sqrt{3})^{5/2}}[231+76\sqrt{3}]+\frac{4}{(2+\sqrt{3})^{5/2}}[231-76\sqrt{3}]\\
&=4(2+\sqrt{3})^{5/2}(231+76\sqrt{3})+4(2-\sqrt{3})^{5/2}(231-76\sqrt{3})\\
&=4(2+\sqrt{3})^{1/2}(2529+1456\sqrt{3})+4(2-\sqrt{3})^{1/2}(2529-1456\sqrt{3})\\
&=4(2+\sqrt{3})^{1/2}[2529+1456\sqrt{3}+(2-\sqrt{3})(2529-1456\sqrt{3})]\\
&=4(2+\sqrt{3})^{1/2}(11,955-3985\sqrt{3})
\endaligned\tag5.5
$$
Substituting (5.5) into the $L$-series (5.4) one obtains
$$
\aligned
L(5,\chi_{24})&=\frac{\pi^5}{2\cdot4!\cdot 12^5}\left[4(2+\sqrt{3})^{1/2}(11,955-3985\sqrt{3})+1225+\frac{68\sqrt{3}}{3} +114\sqrt{2}\right]\\
&=\frac{\pi^5}{2\cdot4!\cdot 12^5\cdot3}\left[12(2+\sqrt{3})^{1/2}(11,955-3985\sqrt{3})+3675+68\sqrt{3}+342\sqrt{2}\right]\\
\endaligned
$$
\flushpar $\therefore\, L(5,\chi_{24})=\frac{\pi^5}{35,831,808}\left[(2+\sqrt{3})^{1/2}(143,460-47,820\sqrt{3})+3675+68\sqrt{3}+342\sqrt{2}\right].$\qed


\Refs\nofrills{References}
\ref \key 1 \by T. J. I'{\smc A} Bromwich\book An introduction to the theory of Infinite Series (3rd edition)\publ Chelsea Publ. Co. \yr (1991)\endref



\ref \key 2 \by K. Kato, N. Kurokawa, T. Saito \book Number Theory I. Fermat's Dream, {\rm Translations of Mathematical Monographs, Vol 186}
\publ Amer. Math. Society \yr (2000) \endref

\ref \key 3 \by H. L. Montgomery, R. C. Vaughan\book Multiplicative number theory I. Classical theory \publ Cambridge Univ. Press \yr (2007)\endref

\ref \key 4 \by G. Sansone, J. Gerretsen \book Lectures on the theory of functions of a complex variable I. Holomorphic Functions
\publ P. Noordhoff Ltd., Groningen, the Netherlands \yr (1960)\endref



\endRefs
\enddocument